\documentclass[11pt]{amsart}
\usepackage[T1]{fontenc}
\usepackage[utf8]{inputenc}
\usepackage[british]{babel}
\usepackage[margin=1in]{geometry}
\usepackage{amsmath,amssymb,amsthm}
\usepackage{palatino,euler,microtype}

\usepackage{booktabs,longtable,array}
\usepackage{graphicx}
\usepackage{tikz}
\usetikzlibrary{positioning,fit,backgrounds,calc,shapes.geometric}
\DeclareMathVersion{jfefpalatino}
\SetSymbolFont{operators}{jfefpalatino}{OT1}{zplm}{m}{n}
\SetSymbolFont{letters}{jfefpalatino}{OML}{zplm}{m}{it}
\SetSymbolFont{symbols}{jfefpalatino}{OMS}{zplm}{m}{n}
\SetSymbolFont{largesymbols}{jfefpalatino}{OMX}{zplm}{m}{n}
\SetMathAlphabet{\mathbf}{jfefpalatino}{OT1}{zplm}{b}{n}
\SetMathAlphabet{\mathit}{jfefpalatino}{OT1}{ppl}{m}{it}
\tikzset{
  jfef text family/.initial=ppl,
  jfef math version/.initial=jfefpalatino,
  jfef font size/.initial=8,
  jfef baseline/.initial=9,
  jfef label font/.style={font={
    \fontfamily{\pgfkeysvalueof{/tikz/jfef text family}}\selectfont
    \fontsize{\pgfkeysvalueof{/tikz/jfef font size}}{\pgfkeysvalueof{/tikz/jfef baseline}}\selectfont
    \mathversion{\pgfkeysvalueof{/tikz/jfef math version}}
  }},
  jfef vertex/.style={
    circle,draw=black!70,line width=.4pt,
    minimum size=5.8mm,inner sep=.8pt,outer sep=0pt,
    fill=white,opacity=1,fill opacity=1,draw opacity=1,text opacity=1,
    text=black,jfef label font
  },
  jfef type/.style={jfef vertex,minimum size=6mm,inner sep=1pt},
  jfef centre/.style={jfef vertex,minimum size=8mm,inner sep=1pt,fill=orange!20},
  jfef auxiliary/.style={jfef vertex,rectangle,minimum size=8mm,inner sep=1pt,fill=cyan!15},
  jfef edge/.style={draw=black!27,line width=.25pt,line cap=round},
  jfef complement edge/.style={draw=black!48,line width=.35pt,line cap=round},
  jfef new edge/.style={draw=blue!55!black,opacity=.55,line width=.35pt,line cap=round},
  jfef new complement edge/.style={draw=blue!65!black,line width=.6pt,line cap=round}
}

\usepackage{xcolor}
\usepackage{float,placeins}
\usepackage{needspace}
\definecolor{darkblue}{rgb}{0,0,0.7}
\definecolor{darkred}{rgb}{0.7,0,0}
\newcommand{\defi}[1]{\textit{\color{blue}#1}}
\usepackage[pagebackref=true]{hyperref}
\hypersetup{colorlinks=true,linkcolor=darkred,citecolor=darkblue,urlcolor=darkblue,
pdftitle={Linear-size graphs at every strict dismantling level}}
\newcommand{\D}{\mathcal D}
\newcommand{\Cl}{\operatorname{Cl}}
\newcommand{\card}[1]{\lvert#1\rvert}

\newcommand{\Nb}[2]{#1[N_{#1}(#2)]}
\newtheorem{theorem}{Theorem}[section]
\newtheorem{lemma}[theorem]{Lemma}
\newtheorem{proposition}[theorem]{Proposition}

\newtheorem{definition}[theorem]{Definition}
\newtheorem{remark}[theorem]{Remark}

\title[Linear-size graphs at every strict dismantling level]{Linear-size graphs at every strict dismantling level}
\author[J. F. Espinoza]{Jes\'us F. Espinoza\textsuperscript{*}}
\thanks{\textsuperscript{*}Corresponding author}
\author[M. E. Fr\'ias-Armenta]{Mart\'in Eduardo Fr\'ias-Armenta}
\author[R. G. Hern\'andez-Amador]{Rosal\'ia G. Hern\'andez-Amador}
\author[H. A. Hern\'andez-Hern\'andez]{H\'ector A. Hern\'andez-Hern\'andez}
\address{Departamento de Matem\'aticas, Universidad de Sonora, Hermosillo, Sonora, M\'exico}
\email[Jes\'us F. Espinoza]{jesusfrancisco.espinoza@unison.mx}
\urladdr[Jes\'us F. Espinoza]{https://orcid.org/0000-0001-5601-3634}
\email[Mart\'in Eduardo Fr\'ias-Armenta]{eduardo.frias@unison.mx}
\urladdr[Mart\'in Eduardo Fr\'ias-Armenta]{https://orcid.org/0000-0002-7687-1270}
\email[Rosal\'ia G. Hern\'andez-Amador]{rosalia.hernandez@unison.mx}
\urladdr[Rosal\'ia G. Hern\'andez-Amador]{https://orcid.org/0000-0001-9657-7012}
\email[H\'ector A. Hern\'andez-Hern\'andez]{hector.hernandez@unison.mx}
\urladdr[H\'ector A. Hern\'andez-Hern\'andez]{https://orcid.org/0000-0002-4046-4841}
\date{}
\begin{document}
\begin{abstract}
Dismantling reduces a graph to one vertex by successive deletions.
At level zero, a vertex is eligible if one of its neighbours is
adjacent to all its other neighbours. At level $k\geq1$, its neighbours
must induce a graph dismantlable at level $k-1$. A graph has strict
level $k$ if it admits a reduction at level $k$ and none at level $k-1$.
For every $k\geq1$, we construct a graph of strict level $k$ with
$3k+8$ vertices. Each step adds three vertices: the neighbours of one
induce the preceding graph, and the other two provide the adjacencies
used to rule out other first deletions at the lower level.
A second construction gives a graph of order $4k+4$ at strict level
$k$. Together the two families give the upper bound
$\min\{4k+4,3k+8\}$ on the least order at strict level $k$.
\end{abstract}
\maketitle
\section{Introduction}
In an ordinary dismantling of a graph, we delete vertices one at a
time until only one remains. A deletion is allowed if one of the
vertex's neighbours is adjacent to all its other neighbours in the
current graph. Such a deletion is called a \defi{fold}. Fieux and Jouve
take ordinary dismantling as level zero and define the higher levels
recursively~\cite{FJ2020}. At level
$k\geq1$, a vertex may be deleted when the subgraph induced by its
neighbours can itself be dismantled at level $k-1$. We call a graph
\defi{strictly $k$-dismantlable} if it can be dismantled at level $k$ and
cannot be dismantled at level $k-1$.

How many vertices are needed for a graph to have strict level $k$?
Fieux and Jouve showed that the levels are distinct by constructing
graphs called cubions. Their example at strict level $k$ has
$2^{k+1}+2k+2$ vertices for every $k\geq1$
\cite{FJ2020}.
The order of these examples grows exponentially with $k$.
We show that linear growth suffices: for every $k\geq1$, we construct
a strictly $k$-dismantlable graph $F_k$ on $3k+8$ vertices.
A second family gives graphs on $4k+4$ vertices, which are smaller
than $F_k$ at the first three levels. Together the families give the upper bound
$\min\{4k+4,3k+8\}$ on the least order at strict level $k$.

Constructing a graph at the next strict level requires two things.
We must exhibit a reduction at the new level, and we must rule out
every possible first deletion at the preceding level.
To start the reduction, we add a vertex, called the \defi{centre},
whose neighbours induce the preceding graph. If that graph has
strict level $k$, we may therefore delete the centre at level $k+1$.
At level $k$ the deletion is forbidden, since the neighbourhood
graph cannot be dismantled at level $k-1$.
After deleting the centre we still have to dismantle the remaining
graph. To exclude level $k$, we also need to examine every other
vertex as a possible first deletion.

For the exclusion, we count cliques in the graph induced by a
vertex's neighbours. These cliques are the faces of its \defi{clique complex}.
Every graph in the dismantling hierarchy has a clique complex of
Euler characteristic one. Applied to the neighbourhood graph, this
gives a test: if the characteristic differs from one, the vertex is
ineligible at every level. We choose the adjacencies of
the two other new vertices so that this obstruction holds away from
the centre in the first layer. A recurrence shows that the relevant
differences from one alternate in sign as layers are added. We also
give explicit folds for the graph left after deleting a centre when
there are at least two layers. These two arguments reduce the proof
for all layers to checks on the starting graph and its first
extension (Theorem~\ref{thm:template}).

The Euler obstruction has a topological explanation. Ordinary folds
correspond to strong collapses of clique complexes, and a graph
belongs to some dismantling class exactly when its clique complex is
non-evasive~\cite{FJ2020}; see also~\cite{BM2012,FL2012}.
A non-evasive complex is collapsible~\cite{KSS1984}.

We begin with the first layer of the $3k+8$ family to show how the
three new vertices are attached. We then give the general criterion
and verify it for the eleven-vertex starting graph. Any construction
satisfying the hypotheses of this Euler criterion must add at least
three vertices per layer. The same folds and clique counts also
determine the effect of deleting a single vertex from $F_k$ for
$k\geq2$. Deleting a centre leaves an ordinarily dismantlable graph.
Deleting any other vertex leaves a graph outside all dismantling
classes. We close with the second family and compare the orders of
the constructions.

\section{Preliminaries}\label{sec:prelim}
All graphs are finite, simple and undirected. For a graph \(H\), write
\(V(H)\) and \(E(H)\) for its vertex and edge sets and \(H[S]\) for the
subgraph induced by \(S\subseteq V(H)\). In particular,
\(H-v=H[V(H)\setminus\{v\}]\). The \defi{open and closed neighbourhood
sets} of \(v\) are
\[
 N_H(v)=\{u\in V(H):uv\in E(H)\},\qquad N_H[v]=N_H(v)\cup\{v\}.
\]
The \defi{open neighbourhood graph} \(H[N_H(v)]\) is the graph
tested by the deletion rule. We reserve the word \emph{link} for
simplicial complexes.

A \defi{cone} is a non-empty graph with a universal vertex, called an
\defi{apex}. The one-vertex graph \(K_1\) is a cone, and the empty graph
is not. A vertex \(u\) is \defi{dominated} by \(t\ne u\) if
\(N_H[u]\subseteq N_H[t]\). Deleting \(u\) under this condition is a
\defi{fold}, written \(u\to t\). We can recognise a fold either
from these closed neighbourhoods or by finding an apex in the
open neighbourhood graph.

\begin{lemma}\label{lem:dom}
The graph \(H[N_H(u)]\) is a cone if and only if \(u\) is dominated
by another vertex.
\end{lemma}
\begin{proof}
Suppose \(N_H[u]\subseteq N_H[t]\) with \(t\ne u\). Since \(u\)
belongs to \(N_H[t]\), the vertices \(u\) and \(t\) are adjacent.
The same inclusion makes \(t\) adjacent to every vertex of
\(N_H(u)\setminus\{t\}\), so \(t\) is an apex of \(H[N_H(u)]\).
Conversely, if \(t\) is an apex of \(H[N_H(u)]\), it is adjacent
to \(u\) and to every other vertex of \(N_H(u)\). Together with
\(t\in N_H[t]\), this gives \(N_H[u]\subseteq N_H[t]\).
For an isolated vertex the open neighbourhood graph is empty, and
there is no vertex that dominates it.
\end{proof}

The cone test supplies the base of the recursive definition.
At each higher level the same kind of test is applied to the graph
on the neighbours, using the class from the preceding level.

\begin{definition}\label{def:hierarchy}
Let \(\D_{-1}\) be the class of non-empty cones. For each integer
\(j\geq0\), a non-empty graph \(H\) belongs to \(\D_j\) if it can
be reduced to \(K_1\) by successive vertex deletions such that the
open neighbourhood graph of each deleted vertex, in the graph remaining
at that step, belongs to \(\D_{j-1}\). A vertex whose open neighbourhood
graph lies in \(\D_{j-1}\) is \defi{\(j\)-dismantlable}.
A graph in \(\D_j\setminus\D_{j-1}\) is \defi{strictly
\(j\)-dismantlable}. The empty graph belongs to none of these classes.
\end{definition}

This is the hierarchy of Fieux and Jouve
\cite{FJ2020}, also used in
\cite{DEFH2023}.
By Lemma~\ref{lem:dom}, the deletions allowed in \(\D_0\) are folds,
so \(\D_0\) is the usual class of dismantlable graphs. At level one,
the allowed deletions are the \(s\)-dismantling moves~\cite{BFJ2010}.
The graph \(K_1\) belongs to every \(\D_j\) because its reduction is empty.
Membership is invariant under graph isomorphism. An isomorphism
takes the apex of a cone to an apex and identifies the open
neighbourhood graphs at each step of a reduction. Induction on
\(j\) therefore preserves the admissibility of each deletion.

To compare two successive levels, we test the same deletion
sequence at both levels.

\begin{lemma}\label{lem:inclusions}
The classes satisfy
\[
 \D_{-1}\subseteq\D_0\subseteq\D_1\subseteq\D_2\subseteq\cdots .
\]
\end{lemma}
\begin{proof}
Start with a cone of apex \(z\). Each other vertex can be folded onto
\(z\), which remains an apex as those vertices are deleted. This
gives \(\D_{-1}\subseteq\D_0\).
Now suppose \(\D_{j-1}\subseteq\D_j\) and \(H\in\D_j\).
Every neighbourhood graph tested in a level-\(j\) reduction of
\(H\) belongs to \(\D_{j-1}\) and hence to \(\D_j\).
The same sequence is therefore a level-\((j+1)\) reduction, completing
the induction.
\end{proof}

For a graph with more than one vertex, a reduction must start with
an admissible deletion. To use such a deletion to prove membership,
we must also be able to reduce the graph it leaves.

\begin{lemma}\label{lem:first}
Let \(\card{V(H)}>1\) and \(j\geq0\).
\begin{enumerate}
\item If \(H[N_H(v)]\notin\D_{j-1}\) for every \(v\in V(H)\), then
\(H\notin\D_j\).
\item If \(H[N_H(v)]\in\D_{j-1}\) and \(H-v\in\D_j\) for some
\(v\in V(H)\), then \(H\in\D_j\).
\end{enumerate}
\end{lemma}
\begin{proof}
In the first case the neighbourhood test fails at every vertex,
so no reduction can start. In the second, the test permits deleting
\(v\). The reduction of \(H-v\) to \(K_1\) at the same level then
completes the sequence.
\end{proof}

\section{One layer of the construction}\label{sec:example}
Let $A$ be the graph on $\{0,1,\ldots,10\}$ given in
Table~\ref{tab:seed-adjacency}.\footnote{The seed can also be written
as \(A=H+\{1,2\}\), where \(H=(V(A),E(A)\setminus\{\{1,2\}\})\).}
We will verify in
Section~\ref{sec:seed3k8} that $A$ has strict level one. To extend it
to a graph $B$ of strict level two, we want a vertex that can be deleted
at level two but is ineligible at level one. Add a vertex $p$, the
centre, whose neighbours are exactly the vertices of $A$. Then
$B[N_B(p)]=A$, so $p$ is eligible at level two and ineligible at
level one.

Add two more vertices $x,y$ to control the neighbourhoods of the
other vertices. Join $x$ to every vertex of $A$ except $2$, join $y$
to every vertex of $A$ except $1$, and add the edge $xy$. Neither
$x$ nor $y$ is adjacent to $p$, so its neighbourhood stays the same.
The resulting graph $B$ has fourteen vertices and sixty-five edges.
The original graph $A$ has eleven vertices and thirty-three edges.

\begin{figure}[htbp]
\centering\begin{tikzpicture}[x=1cm,y=1cm]
\coordinate (w) at (0,0);
\coordinate (p) at (0,2);
\coordinate (x) at (-3,.6);
\coordinate (y) at (3,.6);
\coordinate (a) at (-3,-1.6);
\coordinate (b) at (3,-1.6);
\draw[jfef new edge] (p)--(w);
\draw[jfef new edge] (x)--(w);
\draw[jfef new edge] (y)--(w);
\draw[jfef new edge] (x)--(a);
\draw[jfef new edge] (y)--(b);
\draw[jfef new edge] (x) .. controls (-1.7,-1.7) and (1.7,-1.7) .. (y);
\draw[jfef new edge] (p) .. controls (-5,2.2) and (-5,-1.2) .. (a);
\draw[jfef new edge] (p) .. controls (5,2.2) and (5,-1.2) .. (b);
\node[ellipse,draw=black!70,fill=white,minimum width=3cm,minimum height=1.3cm,jfef label font] at (w) {$W$};
\node[jfef centre] at (p) {$p$};
\node[jfef auxiliary] at (x) {$x$};
\node[jfef auxiliary] at (y) {$y$};
\node[jfef vertex] at (a) {$1$};
\node[jfef vertex] at (b) {$2$};
\end{tikzpicture}
\caption{Edges incident with the three new vertices in the first layer.
The ellipse $W$ represents the nine vertices $V(A)\setminus\{1,2\}$.
An edge to $W$ represents adjacency to all nine. The edges within $A$,
including $\{1,2\}$, are given in Table~\ref{tab:seed-adjacency}.}
\label{fig:first-layer}
\end{figure}
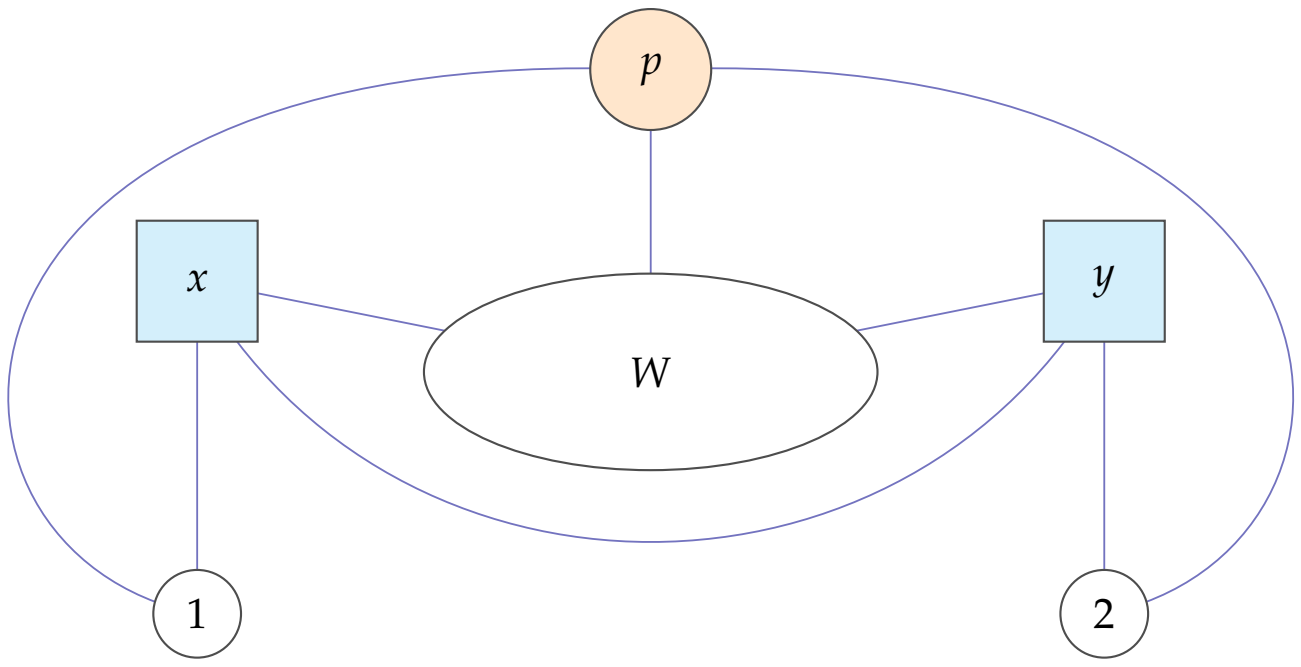

Deleting $p$ is only the first step of a reduction of $B$.
We will show that $B-p$ can already be dismantled at level one, so
the reduction can be completed at level two. To exclude level one,
we must also rule out a first deletion of any vertex other than $p$.
For each such vertex, we will check that the clique complex of its
neighbourhood graph has Euler characteristic different from one.

\begin{table}[htbp]
\centering
\caption{Complete adjacency data for the seed $A$.}\label{tab:seed-adjacency}
\begin{tabular}{cl}
\toprule
\(v\)&\(N_A(v)\)\\\midrule
$0$ & $\{3,4,6,7,8\}$\\
$1$ & $\{2,4,5,6,7,9,10\}$\\
$2$ & $\{1,5,6,7,8,9,10\}$\\
$3$ & $\{0,6,7,8,9,10\}$\\
$4$ & $\{0,1,7,8,9,10\}$\\
$5$ & $\{1,2,6,8,9\}$\\
$6$ & $\{0,1,2,3,5,10\}$\\
$7$ & $\{0,1,2,3,4,9\}$\\
$8$ & $\{0,2,3,4,5,10\}$\\
$9$ & $\{1,2,3,4,5,7\}$\\
$10$ & $\{1,2,3,4,6,8\}$\\
\bottomrule
\end{tabular}

\end{table}

To repeat the construction, we record which vertices of $A$ are
adjacent to each of $x$ and $y$. These two sets are the \defi{profiles}:
\begin{equation}\label{eq:seed-profiles}
 P_0=V(A)\setminus\{2\},\qquad
 P_1=V(A)\setminus\{1\}.
\end{equation}
In later layers, the auxiliary vertices attached according to
$P_0$ and $P_1$ will be of types $0$ and $1$, respectively. The
general construction uses $d$ for the number of types and $R$ for
the graph that records adjacency between types. Here $d=2$, and
the edge $xy$ gives $R=K_2$. For the first layer, the centre and its
two auxiliaries are $p=c_0$, $x=x_{0,0}$ and $y=x_{1,0}$.

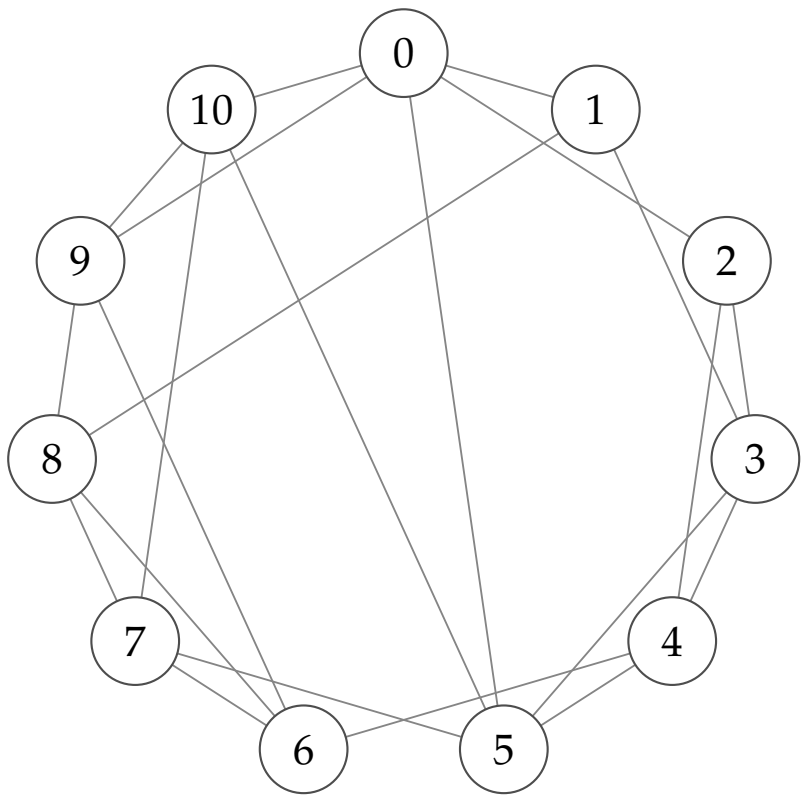
\begin{figure}[htbp]
\centering\begin{tikzpicture}[x=1cm,y=1cm]
\path[use as bounding box] (-3.35,-2.9) rectangle (3.35,2.9);
\coordinate (v0) at (0.000000,2.350000);
\coordinate (v1) at (1.270506,1.976946);
\coordinate (v2) at (2.137635,0.976225);
\coordinate (v3) at (2.326080,-0.334440);
\coordinate (v4) at (1.776011,-1.538923);
\coordinate (v5) at (0.662072,-2.254808);
\coordinate (v6) at (-0.662072,-2.254808);
\coordinate (v7) at (-1.776011,-1.538923);
\coordinate (v8) at (-2.326080,-0.334440);
\coordinate (v9) at (-2.137635,0.976225);
\coordinate (v10) at (-1.270506,1.976946);
\draw[jfef complement edge] (v0) -- (v1);
\draw[jfef complement edge] (v0) -- (v2);
\draw[jfef complement edge] (v0) -- (v5);
\draw[jfef complement edge] (v0) -- (v9);
\draw[jfef complement edge] (v0) -- (v10);
\draw[jfef complement edge] (v1) -- (v3);
\draw[jfef complement edge] (v1) -- (v8);
\draw[jfef complement edge] (v2) -- (v3);
\draw[jfef complement edge] (v2) -- (v4);
\draw[jfef complement edge] (v3) -- (v4);
\draw[jfef complement edge] (v3) -- (v5);
\draw[jfef complement edge] (v4) -- (v5);
\draw[jfef complement edge] (v4) -- (v6);
\draw[jfef complement edge] (v5) -- (v7);
\draw[jfef complement edge] (v5) -- (v10);
\draw[jfef complement edge] (v6) -- (v7);
\draw[jfef complement edge] (v6) -- (v8);
\draw[jfef complement edge] (v6) -- (v9);
\draw[jfef complement edge] (v7) -- (v8);
\draw[jfef complement edge] (v7) -- (v10);
\draw[jfef complement edge] (v8) -- (v9);
\draw[jfef complement edge] (v9) -- (v10);
\node[jfef vertex] at (v0) {$0$};
\node[jfef vertex] at (v1) {$1$};
\node[jfef vertex] at (v2) {$2$};
\node[jfef vertex] at (v3) {$3$};
\node[jfef vertex] at (v4) {$4$};
\node[jfef vertex] at (v5) {$5$};
\node[jfef vertex] at (v6) {$6$};
\node[jfef vertex] at (v7) {$7$};
\node[jfef vertex] at (v8) {$8$};
\node[jfef vertex] at (v9) {$9$};
\node[jfef vertex] at (v10) {$10$};
\end{tikzpicture}
\caption{The complement of $A$, showing the twenty-two pairs that
are not edges of $A$. The vertex labels agree with
Table~\ref{tab:seed-adjacency}.}
\label{fig:seed-complement}
\end{figure}

The next section explains how the checks on $A$ and $B$ give
strict levels for all subsequent layers.

\section{Layers and the increase in dismantling level}\label{sec:construction}

In Section~\ref{sec:example}, the neighbourhood graph of the centre
\(p\) is the entire seed \(A\). We now arrange for the same property
at each step. The neighbourhood graph of a centre in a graph with
\(m\) layers will be isomorphic to the graph with \(m-1\) layers.
The adjacencies of the auxiliaries will let us control the other
neighbourhoods. We first define the construction and then give
conditions under which each layer raises the strict dismantling
level by one.

\subsection{The layer construction}
Start with a non-empty graph \(A\) and fix an integer \(d\geq1\).
Each layer has a centre and one auxiliary of each type in
\(T=\{0,\ldots,d-1\}\). To specify the edges from the auxiliaries to the
seed, choose a subset \(P_j\subseteq V(A)\) for each \(j\in T\).
We call \(P_j\) the \defi{profile} of type \(j\). Its vertices are the
neighbours in \(A\) of every auxiliary of that type. To specify the
edges between auxiliaries of distinct types, choose a graph \(R\)
on \(T\).

Let \(G_m\) denote the graph obtained by adding \(m\) layers.
Index the layers by \(L_m=\{0,\ldots,m-1\}\) for \(m\geq1\), and
put \(L_0=\varnothing\). The vertices of \(G_m\) form the disjoint union
\[
 V(G_m)=V(A)\sqcup\{c_i:i\in L_m\}
                \sqcup\{x_{j,i}:j\in T,\ i\in L_m\}.
\]
Here \(c_i\) is the centre of layer \(i\), and \(x_{j,i}\) is its
auxiliary of type \(j\). The following six rules specify every edge:
\begin{enumerate}
\item The vertices of \(A\) induce the original graph \(A\).
\item The centres form a clique and are adjacent to every vertex of \(A\).
\item For each type \(j\), its auxiliaries \(x_{j,i}\), \(i\in L_m\),
form a clique.
\item For distinct types \(j,t\), an auxiliary of type \(j\) is adjacent
to an auxiliary of type \(t\) exactly when \(jt\in E(R)\), regardless
of their layers.
\item The vertex \(a\in V(A)\) is adjacent to \(x_{j,i}\) exactly when
\(a\in P_j\).
\item The centre \(c_i\) is adjacent to \(x_{j,h}\) exactly when \(i\ne h\).
\end{enumerate}
The sixth rule is what makes a centre omit its own layer from its
neighbourhood graph. For later use, write \(B=G_1\), \(p=c_0\) and
\(x_j=x_{j,0}\). Then \(\Nb{B}{p}=A\).
The seed itself is \(G_0=A\). Adding a layer preserves all existing
adjacencies and adds \(d+1\) vertices, so
\begin{equation}\label{eq:general-order}
 |V(G_m)|=|V(A)|+(d+1)m.
\end{equation}

\begin{figure}[htbp]
\centering\begin{tikzpicture}[x=1cm,y=1cm]
\path[use as bounding box] (-3.1,-2.7) rectangle (3.1,2.7);
\coordinate (v0) at (-2.100000,2.000000);
\coordinate (v1) at (0.000000,2.000000);
\coordinate (v2) at (2.100000,2.000000);
\coordinate (v3) at (-2.100000,-0.000000);
\coordinate (v4) at (-2.100000,-2.000000);
\coordinate (v5) at (0.000000,-0.000000);
\coordinate (v6) at (0.000000,-2.000000);
\coordinate (v7) at (2.100000,-0.000000);
\coordinate (v8) at (2.100000,-2.000000);
\draw[jfef edge] (v0) -- (v1);
\draw[jfef edge] (v0) .. controls (-0.700000,2.900000) and (0.700000,2.900000) .. (v2);
\draw[jfef edge] (v0) -- (v5);
\draw[jfef edge] (v0) -- (v6);
\draw[jfef edge] (v0) -- (v7);
\draw[jfef edge] (v0) .. controls (-0.079310,1.318391) and (1.320690,-0.014943) .. (v8);
\draw[jfef edge] (v1) -- (v2);
\draw[jfef edge] (v1) -- (v3);
\draw[jfef edge] (v1) -- (v4);
\draw[jfef edge] (v1) -- (v7);
\draw[jfef edge] (v1) -- (v8);
\draw[jfef edge] (v2) -- (v3);
\draw[jfef edge] (v2) .. controls (1.320690,0.014943) and (-0.079310,-1.318391) .. (v4);
\draw[jfef edge] (v2) -- (v5);
\draw[jfef edge] (v2) -- (v6);
\draw[jfef edge] (v3) -- (v4);
\draw[jfef edge] (v3) -- (v5);
\draw[jfef edge] (v3) -- (v6);
\draw[jfef edge] (v3) .. controls (-0.700000,0.900000) and (0.700000,0.900000) .. (v7);
\draw[jfef edge] (v3) -- (v8);
\draw[jfef edge] (v4) -- (v5);
\draw[jfef edge] (v4) -- (v6);
\draw[jfef edge] (v4) -- (v7);
\draw[jfef edge] (v4) .. controls (-0.700000,-1.100000) and (0.700000,-1.100000) .. (v8);
\draw[jfef edge] (v5) -- (v6);
\draw[jfef edge] (v5) -- (v7);
\draw[jfef edge] (v5) -- (v8);
\draw[jfef edge] (v6) -- (v7);
\draw[jfef edge] (v6) -- (v8);
\draw[jfef edge] (v7) -- (v8);
\node[jfef centre] at (v0) {$c_{0}$};
\node[jfef centre] at (v1) {$c_{1}$};
\node[jfef centre] at (v2) {$c_{2}$};
\node[jfef auxiliary] at (v3) {$x_{0,0}$};
\node[jfef auxiliary] at (v4) {$x_{1,0}$};
\node[jfef auxiliary] at (v5) {$x_{0,1}$};
\node[jfef auxiliary] at (v6) {$x_{1,1}$};
\node[jfef auxiliary] at (v7) {$x_{0,2}$};
\node[jfef auxiliary] at (v8) {$x_{1,2}$};
\end{tikzpicture}
\caption{The new vertices in three layers. Each column is a layer,
with its centre at the top and one vertex of each auxiliary type
below it. The centres form a clique. The auxiliary vertices also
form a clique in this example because $R=K_2$. Each centre is
adjacent to the auxiliaries in the other two columns. The seed $A$
and its edges to the new vertices are omitted.}
\label{fig:layers}
\end{figure}
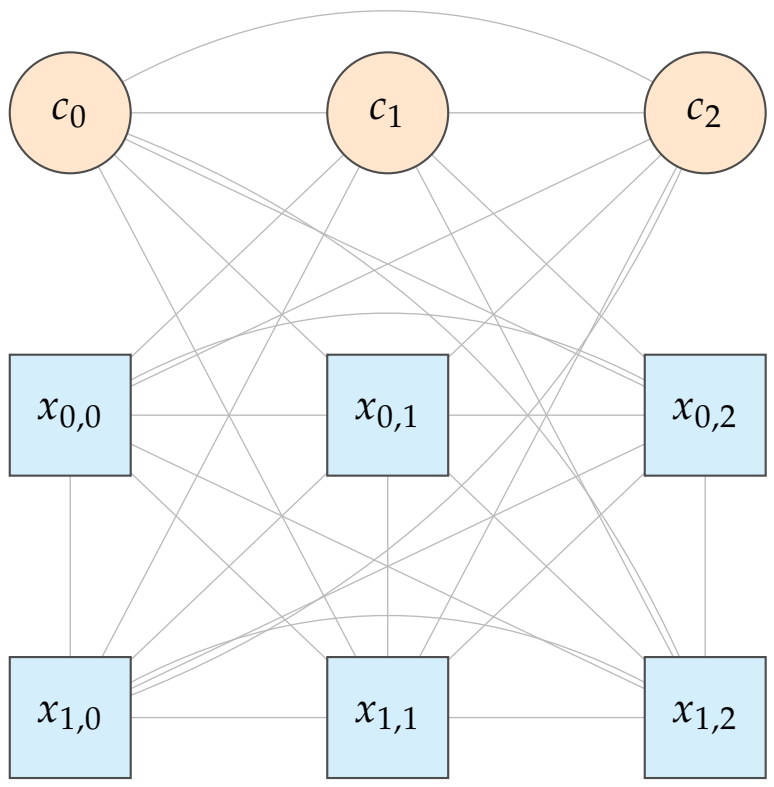

\begin{lemma}\label{lem:centre}
For \(m\geq1\) and \(i\in L_m\), the neighbourhood graph of \(c_i\)
in \(G_m\) is isomorphic to \(G_{m-1}\).
\end{lemma}
\begin{proof}
The neighbours of \(c_i\) are the vertices of \(A\) and the vertices
in layers \(h\in L_m\) with \(h\ne i\):
\[
 N_{G_m}(c_i)=V(A)\sqcup\{c_h:h\ne i\}
                  \sqcup\{x_{j,h}:j\in T,\ h\ne i\}.
\]
Renumber the remaining layers by the increasing bijection
\(\rho_i:L_m\setminus\{i\}\to L_{m-1}\). Fix every vertex of \(A\),
send \(c_h\) to \(c_{\rho_i(h)}\), and send \(x_{j,h}\) to
\(x_{j,\rho_i(h)}\). The first five rules depend on the seed and on
the roles or types of the vertices, which the map preserves. The
sixth depends on whether the two layer indices agree. Since
\(h\ne t\) if and only if \(\rho_i(h)\ne\rho_i(t)\), this rule is
preserved as well. The map is therefore an isomorphism. For \(m=1\)
it is the identity on \(A\).
\end{proof}

To complete a reduction after deleting \(c_i\), we also need to
dismantle the graph that remains. We will fold each auxiliary outside layer \(i\) onto the
auxiliary of the same type in layer \(i\).

\begin{lemma}\label{lem:remainder}
If \(m\geq2\), then \(G_m-c_i\in\D_0\) for every \(i\in L_m\).
\end{lemma}
\begin{proof}
Fix a type \(j\). We will retain \(x_{j,i}\) and fold every
\(x_{j,h}\) with \(h\ne i\) onto it. The two auxiliaries have the
same neighbours in \(A\) and among other types. Within their own
type, their closed neighbourhoods agree because the copies form
a clique. The only difference in their closed neighbourhoods in
\(G_m-c_i\) comes from the centres: \(x_{j,i}\) is adjacent to all
remaining centres, while \(x_{j,h}\) is not adjacent to \(c_h\).
Hence
\[
 N_{G_m-c_i}[x_{j,h}]\subseteq N_{G_m-c_i}[x_{j,i}].
\]
By Lemma~\ref{lem:dom}, this inclusion permits the fold. We perform
these folds for every type. The inclusions persist after earlier
deletions, and all the targets are retained. The vertices left are
those of \(A\), the auxiliaries of layer \(i\), and the \(m-1\) other
centres. Each of these centres is adjacent to every other remaining
vertex. Since \(m\geq2\), there is at least one such centre, so the
remainder is a cone and folds to a single vertex.
\end{proof}

With only one layer, deleting \(p\) leaves no centre to serve as an
apex. We will therefore put a separate hypothesis on \(B-p\).

\subsection{An obstruction that changes sign}
To rule out a non-central vertex as a first deletion, we will check
whether the clique complex of its neighbourhood graph has Euler
characteristic different from one. We express this test through clique
counts. A \defi{clique} is a set of pairwise adjacent vertices. For a graph
\(H\), let \(c_t(H)\) count its cliques of cardinality \(t\), including
the empty clique, so \(c_0(H)=1\). Set
\[
 C_H(z)=\sum_{t\geq0}c_t(H)z^t,
 \qquad q(H)=C_H(-1).
\]
The sum is finite, and \(q(\varnothing)=1\). In terms of the \defi{clique
complex} \(\Cl(H)\), whose faces are the cliques of \(H\),
\begin{equation}\label{eq:euler}
 q(H)=1-\chi(\Cl(H))=-\widetilde\chi(\Cl(H)).
\end{equation}
Here \(\chi\) sums over the non-empty faces: a clique of cardinality
\(t\geq1\) contributes \((-1)^{t-1}\). We use the reduced convention
\(\widetilde\chi=\chi-1\), including
\(\widetilde\chi(\Cl(\varnothing))=-1\).
We will prove that \(q\) vanishes on every dismantling class, so a
non-zero value in a neighbourhood graph rules out deleting its
vertex. The vertex-deletion recurrence for the \defi{clique
polynomial}~\cite{HL1994} lets us follow \(q\) through a reduction.

\begin{lemma}\label{lem:qrec}
For every vertex \(v\) of a graph \(H\),
\[
 C_H(z)=C_{H-v}(z)+zC_{\Nb{H}{v}}(z),
 \qquad q(H)=q(H-v)-q(\Nb{H}{v}).
\]
A cone has \(q=0\), and a fold preserves \(q\).
\end{lemma}
\begin{proof}
Partition the cliques of \(H\) according to whether they contain
\(v\). Those that omit \(v\) are exactly the cliques of \(H-v\).
Removing \(v\) from each clique that contains it gives a clique of
\(\Nb{H}{v}\), and adjoining \(v\) reverses the operation. Since
adjoining \(v\) increases the cardinality by one, these cliques
contribute the term \(zC_{\Nb{H}{v}}(z)\) in the polynomial identity.
Evaluating at \(z=-1\) gives the recurrence for \(q\).

In a cone, pair each clique that omits the apex with the clique
obtained by adjoining the apex. The paired contributions to \(q\)
cancel. By Lemma~\ref{lem:dom}, a vertex \(v\) deleted by a fold has
a cone as its neighbourhood graph. The recurrence thus gives
\(q(H)=q(H-v)\) for that deletion.
\end{proof}

The cone calculation provides the base case for an induction on
the dismantling level.

\begin{lemma}\label{lem:qobstruction}
If \(H\in\D_r\) for an integer \(r\geq-1\), then \(q(H)=0\).
Consequently, \(q(H)\ne0\) excludes \(H\) from every such class.
\end{lemma}
\begin{proof}
The case \(r=-1\) follows from Lemma~\ref{lem:qrec}, since
\(\D_{-1}\) consists of cones. For \(r\geq0\), suppose that \(q\)
vanishes on \(\D_{r-1}\) and take a reduction of \(H\) at level \(r\).
At each deletion, the neighbourhood graph in the current graph
belongs to \(\D_{r-1}\). Its value of \(q\) is zero, so the recurrence
shows that the deletion preserves \(q\). The reduction ends at
\(K_1\), with \(q(K_1)=1-1=0\). Hence \(q(H)=0\).
\end{proof}

A zero value of \(q\) leaves membership undecided. To exclude the
non-central vertices we need non-zero values on their neighbourhood
graphs. It suffices to compute \(q\) on the neighbourhood graphs of
the seed vertices and the auxiliaries in the first extension \(B\),
because adding a layer changes the sign of the corresponding values.
For \(a\in V(A)\) and \(j\in T\), write
\[
 \alpha_a=q(\Nb{B}{a}),\qquad \beta_j=q(\Nb{B}{x_j}).
\]

\begin{lemma}\label{lem:euler-links}
For arbitrary seed, profiles and type graph, and every \(m\geq1\),
\begin{align}
 q(\Nb{G_m}{a})&=(-1)^{m-1}\alpha_a
       &&(a\in V(A)),\label{eq:fixed-euler}\\
 q(\Nb{G_m}{x_{j,i}})&=(-1)^{m-1}\beta_j
       &&(j\in T,\ i\in L_m).\label{eq:aux-euler}
\end{align}
\end{lemma}
\begin{proof}
We argue by induction on \(m\). For \(m=1\), both assertions follow
from the definitions of \(\alpha_a\) and \(\beta_j\). Let \(m\geq2\)
and assume the formulas hold at \(m-1\). In each case \(H\) is
the neighbourhood graph of a non-central vertex. We apply the
recurrence inside \(H\) to a centre of \(G_m\) that lies in \(H\).
We choose that centre so that
deleting it leaves a graph with \(q=0\). The remaining term in
the recurrence is then the negative of \(q\) on the centre's
neighbourhood graph in \(H\).

First fix \(a\in V(A)\) and put \(H=\Nb{G_m}{a}\). Its vertices
are the neighbours of \(a\) in \(A\), all the centres, and all
copies of each type \(j\) with \(a\in P_j\). We may choose any
layer \(i\), since \(H\) contains every centre. In \(H-c_i\),
fold \(x_{j,h}\) onto \(x_{j,i}\) for each of the types present and
every \(h\ne i\). The closed-neighbourhood inclusions from the
proof of Lemma~\ref{lem:remainder} still hold upon restriction to
\(H-c_i\), and every target is present and retained. Once these
folds are done, each of the remaining centres is an apex. At
least one remains because \(m\geq2\). If no auxiliary type is
present, the same observation gives an apex without any folds.
Since folds preserve \(q\) and the resulting cone has \(q=0\),
we have \(q(H-c_i)=0\).

To identify the other term of the recurrence, take the common
neighbours of \(a\) and \(c_i\) in \(G_m\). These are the neighbours
of \(c_i\) in \(H\), and also the neighbours of \(a\) in the
neighbourhood graph of \(c_i\) in \(G_m\). The isomorphism of
Lemma~\ref{lem:centre} fixes \(a\), so the graph induced by these
vertices is isomorphic to \(\Nb{G_{m-1}}{a}\). Thus the recurrence
inside \(H\) gives
\[
 q(\Nb{G_m}{a})=-q(\Nb{G_{m-1}}{a}).
\]
Induction proves~\eqref{eq:fixed-euler}.

Now fix a type \(j\) and a layer \(i\), and put
\(H=\Nb{G_m}{x_{j,i}}\). The vertices of \(H\) are those of
\(A[P_j]\), all copies of the types adjacent to \(j\) in \(R\),
the copies \(x_{j,t}\) with \(t\ne i\), and the centres \(c_t\)
with \(t\ne i\). Choose \(h\in L_m\) with \(h\ne i\). Both \(c_h\)
and \(x_{j,h}\) belong to \(H\). In this list, \(x_{j,h}\) is
adjacent to every seed vertex and every other auxiliary, and
among the centres it omits only \(c_h\). Hence \(x_{j,h}\) is an
apex of \(H-c_h\), and \(q(H-c_h)=0\).

The graph \(\Nb{H}{c_h}\) is induced by the common neighbours of
\(x_{j,i}\) and \(c_h\). In the neighbourhood graph of \(c_h\) in
\(G_m\), these are the neighbours of \(x_{j,i}\).
The isomorphism of Lemma~\ref{lem:centre} renumbers layer \(i\)
as \(\rho_h(i)\), so it identifies this graph with
\(\Nb{G_{m-1}}{x_{j,\rho_h(i)}}\). The recurrence gives
\[
 q(\Nb{G_m}{x_{j,i}})
 =-q(\Nb{G_{m-1}}{x_{j,\rho_h(i)}}).
\]
At \(m-1\) the induction hypothesis applies to the auxiliary of
type \(j\) in every layer, including the layer with index \(\rho_h(i)\).
This proves~\eqref{eq:aux-euler}.
\end{proof}

In both cases, deleting a suitable centre from the neighbourhood
graph leaves a remainder with \(q=0\). The other term of the
recurrence is the value on the corresponding neighbourhood one
layer earlier, taken with a minus sign. In the auxiliary case the
remainder is already a cone, while for a seed vertex folds may be
needed.

\subsection{Strict levels in the layer construction}
The seed determines the level at which we can first delete \(p\).
To finish the reduction of \(B\) we need a hypothesis on \(B-p\),
and to exclude first deletions away from \(p\) we need the
Euler test in \(B\).\footnote{A homological variant of the criterion
uses non-zero reduced homology over $\mathbb F_2$ in the neighbourhood
graphs of non-central vertices. The constructions below use the Euler
criterion.} Lemmas~\ref{lem:remainder}
and~\ref{lem:euler-links} carry the corresponding arguments through
all subsequent layers.

\begin{theorem}\label{thm:template}
Let \(G_m\) be the layer construction above, with \(B=G_1\) and
\(p=c_0\). Suppose that, for some integer \(s\geq0\),
\begin{enumerate}
\item \(A\in\D_s\setminus\D_{s-1}\);
\item \(B-p\in\D_s\);
\item \(q(\Nb{B}{v})\ne0\) for every \(v\in V(B)\setminus\{p\}\).
\end{enumerate}
Then, for every integer \(m\geq0\),
\[
 G_m\in\D_{s+m}\setminus\D_{s+m-1},\qquad
 |V(G_m)|=|V(A)|+(d+1)m.
\]
\end{theorem}
\begin{proof}
We first rule out non-central vertices as candidates for a first
deletion. By hypothesis~(3) and Lemma~\ref{lem:euler-links}, their
neighbourhood graphs in \(G_m\) have non-zero \(q\) whenever
\(m\geq1\). By Lemma~\ref{lem:qobstruction}, none of these graphs
belongs to a dismantling class. Only centres can start a reduction.

For \(m=0\), hypothesis~(1) gives the strict level of the seed.
In \(G_1=B\), the neighbourhood graph of \(p\) is \(A\in\D_s\),
so \(p\) is eligible at level \(s+1\). By hypothesis~(2), the
remainder \(B-p\) belongs to \(\D_s\subseteq\D_{s+1}\), so the
reduction can be completed at level \(s+1\). At level \(s\), the only
candidate for a first deletion is \(p\). Its neighbourhood graph
is outside \(\D_{s-1}\) by hypothesis~(1), so the deletion is
forbidden. For \(s=0\) this says that the neighbourhood graph is
not a non-empty cone. Since \(B\) has more than one vertex,
Lemma~\ref{lem:first} gives \(B\notin\D_s\).

Now suppose \(m\geq2\) and assume that \(G_{m-1}\) has strict level
\(s+m-1\). By Lemma~\ref{lem:centre}, the neighbourhood graph of
any centre of \(G_m\) is isomorphic to \(G_{m-1}\). We may therefore
delete a centre at level \(s+m\) and then use the folds of
Lemma~\ref{lem:remainder} to complete the reduction, since
\(\D_0\subseteq\D_{s+m}\). At level \(s+m-1\), a centre would
need its neighbourhood graph in \(\D_{s+m-2}\). The inductive
hypothesis rules this out, and the Euler obstruction has already
excluded the non-central vertices. Thus there is no eligible first
deletion, and \(G_m\) has more than one vertex.
Lemma~\ref{lem:first} gives \(G_m\notin\D_{s+m-1}\).
The order formula is~\eqref{eq:general-order}.
\end{proof}

\begin{remark}\label{rem:rest-hypothesis}
The proof uses the remainder hypothesis only to obtain
\(B-p\in\D_{s+1}\) in the first step. Thus the conclusion also
holds if hypothesis~(2) is replaced by this weaker condition, while
hypotheses~(1) and~(3) are unchanged.
\end{remark}

Each layer contributes a centre and \(d\) auxiliaries. The Euler
criterion itself puts a lower bound on \(d\). If the seed and
the remainder both belong to dismantling classes, the neighbourhood
graph of a sole auxiliary has \(q=0\).

\begin{proposition}\label{thm:cost-per-layer}
Any data satisfying the hypotheses of Theorem~\ref{thm:template}
have \(d\geq2\). Thus this construction with its Euler obstruction
adds at least three vertices per layer. The same conclusion holds
under the weaker remainder hypothesis of Remark~\ref{rem:rest-hypothesis}.
\end{proposition}
\begin{proof}
The definition takes \(d\geq1\). Even if we allowed \(d=0\),
the first extension would consist of \(A\) and its new universal
centre \(p\). Since \(A\) is non-empty, choose \(a\in V(A)\).
Then \(p\) is an apex of \(\Nb{B}{a}\), including when \(a\)
has no neighbours in \(A\). By Lemma~\ref{lem:qrec},
\(q(\Nb{B}{a})=0\), in contradiction to hypothesis~(3).

For \(d=1\), the sole auxiliary \(x_0\) has neighbourhood graph
\(A[P_0]\) both in \(B\) and in \(B-p\). Deleting \(x_0\) from
\(B-p\) leaves \(A\). Applying Lemma~\ref{lem:qrec} to this
deletion gives
\[
 q(B-p)=q(A)-q(A[P_0]).
\]
By hypotheses~(1) and~(2), both \(A\) and \(B-p\) belong to
dismantling classes, so \(q(A)=q(B-p)=0\) by
Lemma~\ref{lem:qobstruction}. The displayed identity then gives
\(q(\Nb{B}{x_0})=q(A[P_0])=0\), in contradiction to hypothesis~(3).
The same calculation applies when \(B-p\in\D_{s+1}\), since \(q\)
vanishes on that class as well.
\end{proof}

This lower bound is on the number of vertices added by the layer
construction under its Euler hypotheses. In the next section we
verify those hypotheses for the example with two auxiliary types.
The resulting family in Theorem~\ref{thm:family3k8} attains the
bound with three vertices per layer.

\section{Verifying the eleven-vertex seed}\label{sec:seed3k8}
To apply Theorem~\ref{thm:template}, we return to the graphs \(A\) and
\(B\) of Section~\ref{sec:example}. We first establish the strict level
of \(A\). We then give a fold reduction of \(B-p\) and check that
the neighbourhood graphs of the vertices other than \(p\) have
non-zero \(q\). All three checks use the specified adjacencies of \(A\) and \(B\).

\Needspace{6\baselineskip}
\subsection{The seed is strictly \texorpdfstring{$1$}{1}-dismantlable}
We will use folds in two different graphs. To show that \(A\) belongs
to level one, we fold each neighbourhood graph needed in a reduction.
To exclude level zero, we examine the folds available in \(A\) itself.

\begin{lemma}\label{lem:seed-in-D1}
The graph \(A\) belongs to \(\D_1\).
\end{lemma}
\begin{proof}
We delete the vertices of \(A\) in the order
\[
 5,2,1,4,0,6,7,8,9,3,
\]
leaving \(10\). Each row of Table~\ref{tab:seed-reduction} justifies
one deletion. The second column lists the neighbours of the vertex to
be deleted, after the preceding deletions from \(A\). The third column
gives the folds in the graph induced by those neighbours, and the last
column gives the surviving vertex.

To check a listed fold \(u\to t\), let \(H\) be the neighbourhood
graph for that row after its preceding folds. The adjacencies in
Table~\ref{tab:seed-adjacency} give \(N_H[u]\subseteq N_H[t]\), so the
fold is valid by Lemma~\ref{lem:dom}. In the final two rows the
neighbourhood graph is already a single vertex. Thus each row gives a
reduction of its neighbourhood graph to one vertex, which puts that
graph in \(\D_0\) by Definition~\ref{def:hierarchy}. The ten deletions
therefore dismantle \(A\) at level one.
\end{proof}

\begin{table}[htbp]
\centering\small\setlength{\tabcolsep}{4pt}
\caption{A complete level-one reduction of \(A\).}\label{tab:seed-reduction}
\begin{tabular}{cclc}
\toprule
Deleted vertex&Current neighbourhood&Folds within that neighbourhood&Survivor\\\midrule
$5$&$\{1,2,6,8,9\}$&$1\to2,\ 6\to2,\ 8\to2,\ 2\to9$&$9$\\
$2$&$\{1,6,7,8,9,10\}$&$6\to1,\ 7\to1,\ 8\to10,\ 9\to1,\ 1\to10$&$10$\\
$1$&$\{4,6,7,9,10\}$&$6\to10,\ 7\to4,\ 9\to4,\ 4\to10$&$10$\\
$4$&$\{0,7,8,9,10\}$&$9\to7,\ 7\to0,\ 0\to8,\ 8\to10$&$10$\\
$0$&$\{3,6,7,8\}$&$6\to3,\ 7\to3,\ 3\to8$&$8$\\
$6$&$\{3,10\}$&$3\to10$&$10$\\
$7$&$\{3,9\}$&$3\to9$&$9$\\
$8$&$\{3,10\}$&$3\to10$&$10$\\
$9$&$\{3\}$&\text{none}&$3$\\
$3$&$\{10\}$&\text{none}&$10$\\
\bottomrule
\end{tabular}

\end{table}
\FloatBarrier

\Needspace{7\baselineskip}
\begin{lemma}\label{lem:seed-not-D0}
The graph \(A\) does not belong to \(\D_0\).
\end{lemma}
\begin{proof}
To determine whether a level-zero reduction can start, we list the
folds available in \(A\). By Lemma~\ref{lem:dom}, a fold of \(u\) onto
\(t\) in the current graph \(H\) is permitted exactly when \(u\ne t\)
and \(N_H[u]\subseteq N_H[t]\). Comparing the closed neighbourhoods from
Table~\ref{tab:seed-adjacency} gives the complete lists below, first
in \(A\) and then after each forced deletion:
\[
\begin{array}{c|c}
 H&\text{all permitted folds}\\\hline
 A&5\to2\\
 A-5&9\to7\\
 A-\{5,9\}&\text{none}.
\end{array}
\]
Any level-zero reduction would have to delete \(5\) and then \(9\).
The remaining graph on \(\{0,1,2,3,4,6,7,8,10\}\) admits no fold,
so the reduction stops before reaching \(K_1\).
\end{proof}

\subsection{A dismantlable remainder}
The neighbourhood graph of \(p\) is \(A\), so \(p\) can be deleted at
level two. Theorem~\ref{thm:template} uses the hypothesis that the
remaining graph \(B-p\) can be dismantled at level one. In fact,
ordinary folds already suffice.

\begin{lemma}\label{lem:seed-rest}
The graph \(B-p\) belongs to \(\D_0\).
\end{lemma}
\begin{proof}
After removing \(p\), the remaining vertices are \(0,\ldots,10,x,y\).
Perform the following twelve folds in order:
\begin{equation}\label{eq:seed-rest-folds}
\begin{aligned}
 &0\to x,\quad3\to x,\quad4\to x,\quad7\to9,\quad8\to y,\quad9\to5,\\
 &1\to6,\quad2\to6,\quad5\to6,\quad6\to10,\quad10\to x,\quad x\to y.
\end{aligned}
\end{equation}
At each step, the adjacency table and the specified adjacencies
of \(x,y\) give the closed-neighbourhood inclusion for the indicated
fold in the graph left by the preceding folds. Only \(y\) survives.
Thus \(B-p\in\D_0\subseteq\D_1\), giving the required hypothesis of
Theorem~\ref{thm:template} for \(s=1\).
\end{proof}

\Needspace{5\baselineskip}
\subsection{Euler obstructions in the first extension}
A first deletion at level one would require a neighbourhood graph
in \(\D_0\). The centre fails this test because its neighbourhood
graph is \(A\). For each other vertex, we use Lemma~\ref{lem:qobstruction}:
a non-zero \(q\) excludes its neighbourhood graph from \(\D_0\).

\begin{lemma}\label{lem:seed-euler}
Every non-central vertex \(v\) of \(B\) satisfies
\(q(B[N_B(v)])\ne0\).
\end{lemma}
\begin{proof}
For each of the thirteen vertices \(v\ne p\), let \(H_v=B[N_B(v)]\).
To evaluate \(q(H_v)\), we count the cliques of \(H_v\) by size.
In Table~\ref{tab:seed-cliques}, \(c_r(H_v)\) counts the cliques of
size \(r\). The empty clique gives \(c_0(H_v)=1\), and none of the
neighbourhood graphs has a clique of size greater than four.

To obtain a row directly from the adjacency data, enumerate the
subsets of \(N_B(v)\) and count those in which every pair of vertices
is adjacent. These counts are the coefficients of \(C_{H_v}\).
Alternatively, compute the same polynomial using the recurrence in
Lemma~\ref{lem:qrec}, starting from \(C_{\varnothing}(z)=1\).
Each step uses graphs with fewer vertices, so the recursion terminates.
In either case, evaluation at \(z=-1\) gives the last column of
Table~\ref{tab:seed-cliques}, where every entry is non-zero.
\end{proof}

\begin{table}[htbp]
\centering
\caption{All clique counts of \(H_v=B[N_B(v)]\) for \(v\ne p\).}\label{tab:seed-cliques}
\begin{tabular}{crrrrrr}
\toprule
\(v\)&\(c_0(H_v)\)&\(c_1(H_v)\)&\(c_2(H_v)\)&\(c_3(H_v)\)&\(c_4(H_v)\)&\(q(H_v)\)\\\midrule
$0$&1&8&21&20&5&$-1$\\
$1$&1&9&25&24&6&$-1$\\
$2$&1&9&24&21&4&$-1$\\
$3$&1&9&25&24&6&$-1$\\
$4$&1&9&25&24&6&$-1$\\
$5$&1&8&20&16&2&$-1$\\
$6$&1&9&24&21&4&$-1$\\
$7$&1&9&25&25&7&$-1$\\
$8$&1&9&24&21&4&$-1$\\
$9$&1&9&25&24&6&$-1$\\
$10$&1&9&25&23&4&$-2$\\
$x$&1&11&35&37&11&$-1$\\
$y$&1&11&35&36&10&$-1$\\
\bottomrule
\end{tabular}

\end{table}

\subsection{The family at every strict level}
The checks on \(A\) and \(B\) now let us apply Theorem~\ref{thm:template}
to any number of layers, including zero layers.

\begin{theorem}\label{thm:family3k8}
For every integer \(k\geq1\), the graph \(F_k=G_{k-1}\)
with \(d=2\) and \(R=K_2\), specified by Table~\ref{tab:seed-adjacency},
\eqref{eq:seed-profiles}, and the adjacency rules of
Section~\ref{sec:construction} satisfies
\[
 F_k\in\D_k\setminus\D_{k-1},\qquad |V(F_k)|=3k+8.
\]
\end{theorem}
\begin{proof}
Lemmas~\ref{lem:seed-in-D1} and~\ref{lem:seed-not-D0} establish
\(A\in\D_1\setminus\D_0\). The second hypothesis follows from
\(B-p\in\D_0\subseteq\D_1\) by Lemma~\ref{lem:seed-rest}, and the
last is exactly Lemma~\ref{lem:seed-euler}.

Apply Theorem~\ref{thm:template} with \(s=1\) and \(d=2\).
It gives strict level \(m+1\) and order \(11+3m\) after \(m\) layers,
for every \(m\geq0\). Taking \(m=k-1\) gives level \(k\) and order
\(11+3(k-1)=3k+8\). The first two graphs are \(F_1=A\) and \(F_2=B\).
\end{proof}

\section{Deleting a vertex from the main family}\label{sec:consequences}
Can a single vertex deletion preserve the strict level of \(F_k\)?
The proof of strictness already provides two ways to examine the
remaining graph. After deleting a centre we have an explicit fold
reduction. For any other vertex, the Euler recurrence relates the value
of \(q\) on the remaining graph to the value on its neighbourhood graph.

\begin{proposition}\label{prop:single-deletion}
Let $k\geq2$ and let $F_k$ be the graph of
Theorem~\ref{thm:family3k8}. Deleting a centre gives a graph in $\D_0$.
Deleting any other vertex gives a graph outside every class $\D_j$,
$j\geq-1$. Consequently no single vertex deletion preserves its
strict level $k$.
\end{proposition}
\begin{proof}
Write \(F_k=G_m\) with \(m=k-1\geq1\). If we delete a centre
when \(m=1\), the remainder is \(B-p\), which folds to one vertex by
Lemma~\ref{lem:seed-rest}. For \(m\geq2\), Lemma~\ref{lem:remainder}
gives a fold reduction after deleting any centre. In both cases
the remainder belongs to \(\D_0\subseteq\D_{k-1}\) and thus cannot
have strict level \(k\).

Now let \(v\) be a vertex other than a centre. By
Lemma~\ref{lem:euler-links}, the non-zero values in
Table~\ref{tab:seed-cliques} remain non-zero as we add layers, so
\(q(F_k[N_{F_k}(v)])\ne0\). The graph \(F_k\) itself belongs to
\(\D_k\) by Theorem~\ref{thm:family3k8}, so
\(q(F_k)=0\) by Lemma~\ref{lem:qobstruction}. Using these two
values in the recurrence of Lemma~\ref{lem:qrec} gives
\[
 q(F_k-v)=q(F_k)+q(F_k[N_{F_k}(v)])
         =q(F_k[N_{F_k}(v)])\ne0.
\]
Lemma~\ref{lem:qobstruction} therefore excludes \(F_k-v\) from
every dismantling class.
\end{proof}

The conclusion concerns one deletion from \(F_k\), starting with
\(F_2=B\). Minimality among all proper induced subgraphs and least
order among all graphs of strict level \(k\) would require separate
arguments.

\section{A second instance and the orders obtained}\label{sec:secondary}

The \(3k+8\) family starts with eleven vertices at strict level one
and adds three vertices at each step. To obtain smaller graphs at
the first three levels, we use a different seed, with sixteen
vertices at strict level three. Two induced subgraphs of this seed
supply levels one and two. Each layer added to the seed consists of
one centre and three auxiliaries.

We specify the seed \(A^{(4)}\) by listing its missing edges. Its
vertex set is \(\{0,\ldots,15\}\), and every pair of distinct vertices
is joined except for the following twenty-three pairs:
\begin{equation}\label{eq:secondary-missing}
\begin{gathered}
 (0,1),(0,4),(0,15),(1,5),\\
 (1,15),(2,4),(2,5),(2,14),\\
 (3,6),(3,7),(3,14),(4,9),\\
 (5,8),(6,7),(6,10),(7,11),\\
 (8,11),(8,13),(9,10),(9,12),\\
 (10,12),(11,13),(12,13).
\end{gathered}
\end{equation}
Take \(d=3\) and \(E(R)=\{\{0,1\},\{0,2\}\}\), so that type \(0\)
is joined to both other auxiliary types. The neighbours in the seed
of an auxiliary vertex of type \(j\) are specified by \(P_j\):
\begin{equation}\label{eq:secondary-profiles}
\begin{aligned}
 P_0&=V(A^{(4)})\setminus\{8,11,14\},\\
 P_1&=V(A^{(4)})\setminus\{11,15\},\\
 P_2&=V(A^{(4)})\setminus\{13\}.
\end{aligned}
\end{equation}
Write \(G_m^{(4)}\) for the graph obtained by adding \(m\) layers to
this seed, and \(B^{(4)}=G_1^{(4)}\) for its first extension. For the
graphs at levels one and two, take the induced subgraphs on the sets
\[
 S_1=\{2,4,5,8,9,10,12,13\},\qquad S_2=\{2,3,\ldots,13\},
\]
respectively. The full family is
\begin{equation}\label{eq:secondary-family}
 H_1=A^{(4)}[S_1],\qquad H_2=A^{(4)}[S_2],\qquad
 H_k=G_{k-3}^{(4)}\quad(k\geq3).
\end{equation}
The proof uses the layer theorem for the graphs starting at the seed
and verifies the two smaller induced subgraphs separately.

\begin{proposition}\label{prop:family4k4}
For every \(k\geq1\),
\[
 H_k\in\D_k\setminus\D_{k-1},\qquad |V(H_k)|=4k+4.
\]
\end{proposition}
\begin{proof}
We first check membership in the dismantling classes for the seed,
the first remainder and the induced subgraphs \(H_1,H_2\).
The empty graph belongs to no dismantling class, every non-empty
cone belongs to all of them, and a non-cone lies outside \(\D_{-1}\).
In the remaining cases a graph \(J\) belongs to \(\D_j\), with
\(j\geq0\), exactly when there is a vertex \(v\) for which
\(J[N_J(v)]\in\D_{j-1}\) and \(J-v\in\D_j\).
Necessity follows by taking the first vertex of a reduction in
Definition~\ref{def:hierarchy}, and sufficiency follows from
Lemma~\ref{lem:first}. Both conditions involve smaller graphs, so
the rule gives a finite recursion. In these recursive cases, a
certificate of membership specifies a vertex \(v\) together with
certificates of both \(J[N_J(v)]\in\D_{j-1}\) and \(J-v\in\D_j\).
Non-membership is certified by establishing recursively, for every
\(v\in V(J)\), that either \(J[N_J(v)]\notin\D_{j-1}\) or
\(J-v\notin\D_j\).

For the graphs specified above, the recursive checks give
\[
 A^{(4)}\in\D_3\setminus\D_2,\qquad B^{(4)}-p\in\D_3,\qquad
 H_1\in\D_1\setminus\D_0,\qquad H_2\in\D_2\setminus\D_1 .
\]
To check the Euler hypothesis, apply Lemma~\ref{lem:qrec}, starting
from \(q(\varnothing)=1\), to each of the nineteen non-central
neighbourhood graphs in \(B^{(4)}\). The resulting values are
\(-1,-2,-3,-4,-7\), with multiplicities \(8,6,3,1,1\), respectively.
They are all non-zero. The supplementary material contains the recursive
certificates and clique counts, together with a program that verifies
them.

With \(s=d=3\), Theorem~\ref{thm:template} gives strict level
\(m+3\) and order \(16+4m\) for \(G_m^{(4)}\). Setting \(m=k-3\)
proves the assertion for \(k\geq3\). For \(k=1,2\), the membership
checks above give the strict levels of \(H_1\) and \(H_2\), whose
vertex sets have eight and twelve vertices, respectively.
\end{proof}

Let \(m_k\) be the least order of a strictly \(k\)-dismantlable
graph. Taking the smaller of the orders of \(F_k\) and \(H_k\) gives
the constructive upper bound
\[
 m_k\leq\min\{4k+4,3k+8\}\qquad(k\geq1),
\]
with \(H_k\) smaller at \(k=1,2,3\) and \(F_k\) smaller for every
\(k\geq5\). The orders agree at \(k=4\).
The layer construction also yields an instance of order \(6k+6\),
whose larger order does not improve this bound.

At strict level \(k\), the cubion \(Q_{k+1}\) of Fieux and Jouve
has order \(2^{k+1}+2k+2\)~\cite{FJ2020}.
At \(k=1\), the cubion order and our bound both equal the known
minimum \(m_1=8\)~\cite{DEFH2023}. For \(k\geq2\), our bound is
smaller than the cubion order, since \(4k+4<2^{k+1}+2k+2\).
At \(k=2\), \(F_2\) and the cubion both have fourteen vertices,
so the improvement uses the certified twelve-vertex graph \(H_2\).

\section*{Acknowledgements}
The authors acknowledge the \'Area de C\'omputo de Alto Rendimiento
(ACARUS) at the Universidad de Sonora for providing access to the
YUCA high-performance computing cluster used in this work.

\section*{Declarations}
\noindent\textbf{Funding.}
This work was supported by the Universidad de Sonora through the
internal research funding call of the Direcci\'on de Apoyo a Docentes,
Investigaci\'on y Posgrado (DADIP), project USO315010378.

\smallskip
\noindent\textbf{Competing interests.}
The authors declare that they have no competing interests.

\smallskip
\noindent\textbf{Author contributions.}
All authors contributed equally to this work.

\smallskip
\noindent\textbf{Declaration of generative AI and AI-assisted technologies.}
OpenAI Codex and Anthropic Claude Science assisted with examining the
constructions and arguments, computational verification, figures, and
drafting and editing. The authors are responsible for the content of
the article.

\smallskip
\noindent\textbf{Data and code.}
The article gives the adjacency data and finite checks for the
\(3k+8\) family. Section~\ref{sec:secondary} specifies the graphs
and finite checks for the \(4k+4\) family. The accompanying
verification files contain the graph inputs, recursive certificates
and a checker for those finite assertions.
These files are available on
\href{https://github.com/Geometric-and-Combinatorial-Structures/strictly-k-dismantlable-graphs}{GitHub}
and archived on Zenodo (v1.0.0,
\href{https://doi.org/10.5281/zenodo.24099942}{\mbox{DOI~10.5281/zenodo.24099942}}).

\bibliographystyle{siam}
\bibliography{refs}
\end{document}